\documentclass{article}%
\usepackage{amsmath}
\usepackage{amsfonts}
\usepackage{amssymb}
\usepackage{graphicx}%
\providecommand{\U}[1]{\protect\rule{.1in}{.1in}}
\newtheorem{theorem}{Theorem}

\newtheorem{corollary}[theorem]{Corollary}

\newtheorem{definition}[theorem]{Definition}

\newtheorem{proposition}[theorem]{Proposition}
\newtheorem{remark}[theorem]{Remark}

\newenvironment{proof}[1][Proof]{\noindent\textbf{#1.} }{\ \rule{0.5em}{0.5em}}
\begin{document}

\title{Klein-Dirac Quadric and Multidimensional Toda Lattice \\via Pseudo-positive Moment Problem}
\author{S. Albeverio\\Instititute of Applied Mathematics and HCM, IZKS, BiBoS\\University of Bonn, Endenicherallee 60, 53115 Bonn\\CERFIM (Locarno)\thinspace\ \\and\thinspace\ \\O. Kounchev\\Institute of Mathematics and Informatics \\Bulgarian Academy of Sciences\\Acad. G. Bonchev St. 8, 1113 Sofia\\and IZKS-University of Bonn}
\maketitle

\textbf{ ABSTRACT.} In $1974$ J\"{u}rgen Moser has shown that the classical
Moment Problem plays a fundamental role for the theory of completely
integrable systems, by proving that the simplest case of the finite Toda
lattice is described exhaustively in its terms. In particular, the Jacobi
matrix defined by the Flaschka variables corresponds to the Schr\"{o}dinger
operator related to the Korteweg-de Vries operator. 

In the present paper we use a recent breakthrough in the area of the
multidimensional Moment Problem in order to develop a multidimensional
generalization of the finite Toda lattice. Our construction is based on the
notion of pseudo-positive measure and the notion of multidimensional
Markov-Stieltjes transform naturally defined on the Klein-Dirac quadric. An
important consequence of the properties of the Markov-Stieltjes transform is a
new method for summation of multidimensional divergent series on the
Klein-Dirac quadric. 



\section{Introduction}

In the present paper we provide a novel \emph{multidimensional} generalization
of the finite Toda lattice which is governed by an operator which generalizes
the Schr\"{o}dinger one-dimensional operator. Our construction is based on a
new notion of multidimensional Markov-Stieltjes transform naturally defined on
the Klein-Dirac quadric. Our research is a direct generalization of the work
of J\"{u}rgen Moser, \cite{moser1}, \cite{MoserNotes}, \cite{moser2}, who
applied the one-dimensional Moment Problem to study exhaustively the
non-periodic finite Toda lattice.

Let us note that there has been a large amount of research in the area of
multidimensional non-linear evolutions aimed at generalization of the Toda
lattice and the KdV equation. These have led some researchers
(\cite{bealsCoifman}) to the conclusion that there is a "...lack of
interesting non-linear evolutions associated to the Laplacian. A major open
question is to find those operators and systems that do have associated,
stable computable evolutions". Accordingly, one of the major objectives of the
present research is to create models for such operators.

The structure of the paper is as follows: In section \ref{Sclassic} we recall
some basic facts from the classical Moment Problem and the finite Toda
lattice. In section \ref{Spseudo} we recall the (multidimensional)
pseudo-positive Moment Problem by using the framework of the Klein-Dirac
quadric. In section \ref{SNevan} we generalize the Nevanlinna class of
analytic functions to holomorphic functions on the Klein-Dirac quadric, and we
prove a generalization of the theorem of Riesz-Nevanlinna for encoding the
Moment Problem by the means of function theory. As an interesting by-product
we obtain a new method for summation of multidimensional divergent series on
the Klein-Dirac quadric. In section \ref{Smulti} we define our generalization
as the \textbf{pseudo-positive Toda lattice }and prove conditions for its
existence, Theorem \ref{Tnice}. In section \ref{SJacobi} we describe a
functional model for a generalization of the one-dimensional discrete
Schr\"{o}dinger (Jacobi matrix) operator. Finally, in section \ref{Sspecial}
we present a special class of pseudo-positive Toda lattices.

\section{The classical finite Toda lattice and the Moment
Problem\label{Sclassic}}

The Toda lattice is a system of unit masses, connected by nonlinear springs
governed by an exponential restoring force. The equations of motion are
derived from a Hamiltonian containing the displacements $x_{j}$ of the $j$-th
mass and the corresponding momentum $y_{j}.$ The Toda lattice is a
finite-dimensional analog to the Korteveg-de Vries partial differential
equations, cf. \cite{flaschka0}, \cite{flaschka}. 

In the seminal papers by J\"{u}rgen Moser \cite{moser1}, \cite{MoserNotes},
\cite{moser2}, the finite Toda lattice was considered with \emph{non-periodic
boundary conditions}, and a complete solution was given in terms of the
classical Moment Problem. We provide a short outline of these results since we
will use them essentially in our multidimensional generalization. 

We introduce the Hamiltonian given by
\begin{equation}
H=\frac{1}{2}\sum_{j=1}^{N}y_{j}^{2}+\sum_{j=1}^{N-1}e^{x_{j}-x_{j+1}%
},\label{Hxy}%
\end{equation}
$x_{j}$ and $y_{j}$ being real variables corresponding to the mass
displacements and their momenta. The corresponding equations of motion for all
times $t\in\mathbb{R}$ are
\begin{subequations}
\begin{align}
x_{j}^{\prime} &  =y_{j},\qquad\qquad\qquad\qquad\quad
j=1,2,...,N\label{TodaLattice}\\
y_{j}^{\prime} &  =e^{x_{j-1}-x_{j}}-e^{x_{j}-x_{j+1}},\qquad
j=2,3,...,N-1\nonumber\\
y_{1}^{\prime} &  =-e^{x_{1}-x_{2}}\nonumber\\
y_{N}^{\prime} &  =e^{x_{N-1}-x_{N}},\nonumber
\end{align}
(where $^{\prime}$ denotes time derivative). Obviously, they are equivalent
to
\end{subequations}
\begin{align*}
x_{j}^{\prime\prime} &  =e^{x_{j-1}-x_{j}}-e^{x_{j}-x_{j+1}},\qquad
j=1,2,...,N\\
x_{0} &  :=-\infty,\quad x_{N+1}:=+\infty,
\end{align*}
so that $e^{x_{0}-x_{1}}=0$ and $e^{x_{N}-x_{N+1}}=0.$ Following H. Flaschka
\cite{flaschka}, one introduces the new variables
\begin{align}
a_{j} &  =\frac{1}{2}e^{\left(  x_{j}-x_{j+1}\right)  /2},\qquad
j=1,2,...,N-1\label{FlaschkaChange}\\
b_{j} &  =-\frac{1}{2}y_{j},\qquad j=1,2,...,N,\label{FlaschkaChange2}\\
a_{0} &  =0,\quad a_{N}=0,\label{FlaschkaChange3}%
\end{align}
which satisfy the system
\begin{align}
a_{j}^{\prime} &  =a_{j}\left(  b_{j+1}-b_{j}\right)  ,\qquad
j=1,2,...,N-1\label{ab1}\\
b_{j}^{\prime} &  =2\left(  a_{j}^{2}-a_{j-1}^{2}\right)  ,\qquad\text{
}j=1,2,...,N\label{ab2}\\
a_{0} &  =0,\quad a_{N}=0.\label{ab3}%
\end{align}
Let us note that the map from the domain $D_{1}:=\left\{  \left(  x,y\right)
\in\mathbb{R}^{2N}\right\}  $ to the domain $D\subset\mathbb{R}^{2N-1}$ given
by
\[
D:=\left\{  a\in\mathbb{R}^{N-1},b\in\mathbb{R}^{N}:a_{j}>0,\quad\text{for
}j=1,...,N-1\right\}
\]
is invertible only on equivalence classes in $D_{1}$ defined by the
equivalence relation
\begin{equation}
\left(  x,y\right)  \sim\left(  \widetilde{x},\widetilde{y}\right)  \iff
x_{j}-\widetilde{x}_{j}\quad\text{is independent of }j,\label{configuration}%
\end{equation}
cf. \cite{moser1}, p. $471$. We call these every such class configuration of
$\left(  x,y\right)  .$   

We have then for the Hamiltonian (\ref{Hxy}) the following representation in
terms of the new variables:
\begin{equation}
H=4\left\{
{\displaystyle\sum_{j=1}^{N-1}}
a_{j}^{2}+\frac{1}{2}%
{\displaystyle\sum_{j=1}^{N}}
b_{j}^{2}\right\}  .\label{Hab}%
\end{equation}
The following result holds (see \cite{moser1}, Lemma in Section 2, and (2.5)).

\begin{proposition}
Let us consider the system (\ref{ab1})--(\ref{ab2}) with the requirement
(\ref{ab3}) replaced by the requirement that $a_{0}\left(  t\right)  $ and
$a_{N}\left(  t\right)  $ belong to $L_{2}\left(  \mathbb{R}\right)  $ as
functions of $t\in\mathbb{R},$ i.e. satisfy
\[%
{\displaystyle\int_{-\infty}^{\infty}}
\left(  a_{0}^{2}\left(  t\right)  +a_{N}^{2}\left(  t\right)  \right)
dt<\infty.
\]
Then for every $j=2,...,N$ we have
\[%
{\displaystyle\int_{-\infty}^{\infty}}
\left(  a_{1}^{2}+a_{j-1}^{2}\right)  dt<\infty.
\]
Moreover, for any solution of (\ref{ab1})-(\ref{ab3}) satisfying $a_{j}\left(
t\right)  >0,$ $j=1,2,...,N-1,$ we have
\[
a_{j}\left(  t\right)  \longrightarrow0\qquad\text{for }t\longrightarrow
\pm\infty\quad\text{and }j=1,2,...,N-1,
\]
and
\[
b_{j}\left(  t\right)  \overset{t\rightarrow\infty}{\longrightarrow}%
b_{j}\left(  \infty\right)
\]
for some $b_{j}\left(  \infty\right)  \in\mathbb{R}.$
\end{proposition}


Flaschka \cite{flaschka0}, \cite{flaschka} has found the following Lax
representation
\[
L^{\prime}=BL-LB,
\]
where we have the $N\times N$ (symmetric) Jacobi matrix
\begin{equation}
L=\left(
\begin{array}
[c]{ccccc}%
b_{1} & a_{1} &  &  & 0\\
a_{1} & b_{2} &  &  & \\
&  &  &  & \\
&  &  & b_{N-1} & a_{N-1}\\
0 &  &  & a_{N-1} & b_{N}%
\end{array}
\right)  \label{L}%
\end{equation}
and
\begin{equation}
B=\left(
\begin{array}
[c]{ccccc}%
0 & a_{1} &  &  & 0\\
-a_{1} & 0 &  &  & \\
&  &  &  & \\
&  &  & 0 & a_{N-1}\\
0 &  &  & -a_{N-1} & 0
\end{array}
\right)  .\label{B}%
\end{equation}

\begin{remark}
The matrix $L$ is the analog to the Schr\"{o}dinger operator appearing  in the
Inverse Scattering method for solving the KdV equation. 
\end{remark}

Moreover, he found an orthogonal matrix $U$ which satisfies
\begin{align*}
U^{\prime} &  =BU\\
U\left(  0\right)   &  =I
\end{align*}
for which one has
\[
\left(  U^{-1}LU\right)  ^{\prime}=0
\]
which implies
\[
U^{-1}LU=L\left(  0\right)  .
\]
The latter equality implies that the eigenvalues $\left(  \lambda_{j}\right)
_{j=1}^{N}$ of the symmetric matrix $L,$ which are real and distinct, do not
depend on $t.$ Thus the eigenvalues (and symmetric functions of them) are
first integrals of the Toda flow%
\[
\left(  a_{j}\left(  0\right)  ,b_{j}\left(  0\right)  \right)
\longrightarrow\left(  a_{j}\left(  t\right)  ,b_{j}\left(  t\right)  \right)
\qquad\text{for }j=1,2,...,N
\]
associated with (\ref{ab1})-(\ref{ab3}). In particular, for every integer
$m\geq0,$ the trace $\operatorname*{tr}\left(  L^{m}\right)  $ is a symmetric
function of the eigenvalues. We see that, in particular,
\begin{equation}
\operatorname*{tr}\left(  L^{2}\right)  =%
{\displaystyle\sum_{j=1}^{N}}
\lambda_{j}^{2}=\frac{1}{2}H=2\left\{
{\displaystyle\sum_{j=1}^{N-1}}
a_{j}^{2}+\frac{1}{2}%
{\displaystyle\sum_{j=1}^{N}}
b_{j}^{2}\right\}  .\label{traceL2}%
\end{equation}
(\cite{MoserNotes}, p. $218$ ).

Further Flaschka uses the resolvent matrix
\[
R\left(  \lambda\right)  =\left(  \lambda I-L\right)  ^{-1}\qquad\text{for
}\lambda\notin\sigma\left(  L\right)  ,
\]
with $\sigma\left(  L\right)  $ denoting the spectrum of $L.$

We define the important function $f\left(  \lambda,t\right)  $ by putting
\[
R_{N,N}\left(  \lambda,t\right)  =\left\langle R\left(  \lambda\right)
e_{N},e_{N}\right\rangle =:f\left(  \lambda,t\right)  ,\qquad\lambda
\notin\sigma\left(  L\right)  ,\ t\in\mathbb{R}.
\]
Here $e_{N}$ denotes the $N-$vector $\left(  0,0,...,0,1\right)  .$ A simple
argument (\cite{moser1}) shows that
\begin{equation}
f\left(  \lambda,t\right)  =\sum_{j=1}^{N}\frac{r_{j}^{2}\left(  t\right)
}{\lambda-\lambda_{j}\left(  t\right)  },\qquad\lambda\notin\sigma\left(
L\right)  ,\ t\in\mathbb{R},\label{flambdat}%
\end{equation}
for some appropriate real functions $r_{j}\left(  t\right)  $ and $\lambda
_{j}\left(  t\right)  .$ The function $f$ satisfies
\[
\lambda f\left(  \lambda,t\right)  \longrightarrow1\qquad\text{for }\left\vert
\lambda\right\vert \longrightarrow\infty
\]
which implies
\[%
{\displaystyle\sum_{j=1}^{N}}
r_{j}^{2}\left(  t\right)  =1,\qquad t\in\mathbb{R}.
\]
These $\lambda_{j},$ $r_{j}$ are the new variables which show the direct
relation to the Moment Problem.

Let us remark that a classical result of Stieltjes shows that the function
$f\left(  \lambda,t\right)  $ has the pointwise convergent continued fraction
representation:
\begin{equation}
f\left(  \lambda,t\right)  =\frac{1}{\lambda-b_{N}\left(  t\right)
-\frac{a_{N-1}^{2}\left(  t\right)  }{\lambda-b_{N-1}-\cdot\cdot\cdot}}%
=\frac{Q_{N}\left(  \lambda,t\right)  }{P_{N}\left(  \lambda,t\right)
},\label{CF}%
\end{equation}
in terms of the variables $a_{j},b_{j}$ given by (\ref{ab1})-(\ref{ab3}), cf.
\cite{Akhiezer}, \cite{nikishinSorokin}, \cite{simon98}, \cite{teschl}. 

The matrix $L$ contains the $3-$term recurrence relations which generate the
polynomials $P_{N}\left(  \lambda,t\right)  $ which are orthogonal with
respect to the measure $\mu.$ Let $Q_{N}\left(  \lambda,t\right)  $ be the
\emph{second kind orthogonal polynomials} which are obtained by setting
\[
Q_{N}\left(  \tau,t\right)  =%
{\displaystyle\int}
\frac{P_{N}\left(  \tau,t\right)  -P_{N}\left(  \lambda,t\right)  }%
{\tau-\lambda}d\mu\left(  \lambda,t\right)  ,
\]
where $d\mu\left(  \lambda,t\right)  $ is the $t$ -dependent measure in
$\lambda$ on $\mathbb{R}$ defined by
\begin{equation}
d\mu\left(  \lambda,t\right)  :=%
{\displaystyle\sum_{j=1}^{N}}
r_{j}^{2}\left(  t\right)  \delta\left(  \lambda-\lambda_{j}\left(  t\right)
\right)  , \label{dmu}%
\end{equation}
with $\delta$ being the Dirac function concentrated at the origin.
$d\mu\left(  \lambda,t\right)  $ is in fact the spectral measure for the
Jacobi matrix $L.$

From (\ref{flambdat}) and (\ref{dmu}) we have
\[
f\left(  \lambda,t\right)  =\sum_{j=1}^{N}\frac{r_{j}^{2}\left(  t\right)
}{\lambda-\lambda_{j}\left(  t\right)  }=%
{\displaystyle\int}
\frac{d\mu\left(  \tau,t\right)  }{\lambda-\tau},\qquad\lambda\notin
\sigma\left(  L\right)  ,\ t\in\mathbb{R}.
\]

The main observation by J. Moser is contained in the following Proposition
(see \cite{moser1}, p. $480$); it expresses the dynamics of the variables
$r_{j}$ and $\lambda_{j}$ as a result of the dynamics of the variables $a_{j}$
and $b_{j}.$

\begin{proposition}
The variables $r_{j}\left(  t\right)  ,$ $\lambda_{j}\left(  t\right)  $
associated with the solutions of the classical finite Toda lattice given by
(\ref{ab1})-(\ref{ab2}) satisfy the linear equations
\[
\lambda_{j}^{\prime}\left(  t\right)  =0,\qquad r_{j}^{\prime}\left(
t\right)  =-\lambda_{j}\left(  t\right)  r_{j}\left(  t\right)  ,\qquad
j=1,2,...,N.
\]
The solution $\left(  \lambda_{j},r_{j}\right)  ,$ $j=1,2,...,N$ is given by%
\begin{align}
\lambda_{j}\left(  t\right)   &  =\lambda_{j}\left(  0\right)  :=\lambda
_{j}\label{TodaSolution}\\
r_{j}^{2}\left(  t\right)   &  =\frac{r_{j}^{2}\left(  0\right)
e^{-2\lambda_{j}t}}{%
{\displaystyle\sum_{j=1}^{N}}
r_{j}^{2}\left(  0\right)  e^{-2\lambda_{j}t}},\qquad t\in\mathbb{R}%
.\label{TodaSolution2}%
\end{align}
\qquad\qquad
\end{proposition}

For the details of the above representations we refer to the excellent
exposition on the Classical Moment Problem in \cite{simon98}, and for the
finite Toda lattice to \cite{simonSzego}, \cite{teschl}.

\section{Klein-Dirac quadric and pseudo-positive Moment Problem\label{Spseudo}%
}

A \textbf{central role} in our further consideration will be played by the
\textbf{Klein-Dirac quadric} which generalizes the complex plane $\mathbb{C}.$
It is defined by
\[
\operatorname{KDQ}:=\left\{  \zeta\theta:\zeta\in\mathbb{C},\ \theta
\in\mathbb{S}^{n-1}\right\}  =\left(  \mathbb{C}\times\mathbb{S}^{n-1}\right)
_{\mathbb{Z}_{2}}.
\]
For $n=3$ the set $\operatorname{KDQ}$ has been introduced by Klein in his
famous treatise. The notation $\zeta\theta$ denotes the identification of the
points $\left(  \zeta,\theta\right)  $ and $\left(  -\zeta,-\theta\right)  $
in $\mathbb{C}\times\mathbb{S}^{n-1}$ which is the antipodal identification
$\mathbb{Z}_{2}$ of points in $\mathbb{C}\times\mathbb{S}^{n-1}$ expressed by
the notation $\left(  \mathbb{C}\times\mathbb{S}^{n-1}\right)  _{\mathbb{Z}%
_{2}}.$ The unit ball in $\operatorname{KDQ}$ is given by
\[
B_{1}:=\left\{  \zeta\theta:\zeta\in\mathbb{C},\ \left\vert \zeta\right\vert
<1,\ \theta\in\mathbb{S}^{n-1}\right\}  ;
\]
respectively we define the ball of radius $R$ as
\[
B_{R}:=R\cdot B_{1}=\left\{  R\zeta\theta:\zeta\theta\in B_{1}\right\}  .
\]
The topological boundary $\partial\left(  B_{1}\right)  $ is the
\textbf{compactified} \textbf{Minkowski space-time}%
\[
\left(  \mathbb{S}^{1}\times\mathbb{S}^{n-1}\right)  _{\mathbb{Z}_{2}%
}=\left\{  \zeta\theta:\zeta\in\mathbb{S}^{1},\ \theta\in\mathbb{S}%
^{n-1}\right\}  ,
\]
where again $\zeta\theta$ denotes the identification of $\left(  \zeta
,\theta\right)  $ and $\left(  -\zeta,-\theta\right)  $ in $\mathbb{S}%
^{1}\times\mathbb{S}^{n-1}.$

\begin{remark}
The Klein-Dirac quadric has appeared in the works of P. Dirac while he
extended the equations of physics from the Minkowski space-time to a
$4$-dimensional conformal space, in order "to bring a greater symmetry of the
equations into evidence", \cite{dirac1}, \cite{dirac2}. Apparently, the term
Klein-Dirac space has been coined in Conformal Field Theory by Ivan Todorov,
cf.  \cite{nikolovTodorov}, and references therein. The Klein-Dirac quadric
appears in a natural way as a space where the solutions of elliptic equations
(polyharmonic functions) are extended and enjoy a generalized Cauchy type
formula, cf. \cite{aron}, \cite{kounchevRenderResultate},
\cite{kounchevRenderPliska}. 
\end{remark}

In the previous section we saw the key role of the function
\begin{equation}
f\left(  \lambda\right)  =%
{\displaystyle\int}
\frac{d\mu\left(  \tau\right)  }{\lambda-\tau}\qquad\text{for }\lambda
\notin\sigma\left(  L\right)  \label{Stieltjes-transform}%
\end{equation}
which is the Stieltjes transform (sometimes also called Markov transform, cf.
\cite{nikishinSorokin}) of the measure $\mu$ defined by (\ref{dmu}) and is
obtained from the resolvent $R\left(  \lambda\right)  =\left(  \lambda
I-L\right)  ^{-1}$ of the Jacobi operator $L,$ $\lambda\notin\sigma\left(
L\right)  .$

A main achievement of the multidimensional pseudo-positive Moment Problem
introduced in \cite{kounchevrenderARXIV}, \cite{kounchevrender} is that it has
at its disposal analogs to all notions available in the one-dimensional case.
In particular, we have a generalization of the Stieltjes transform (called
here Stieltjes-Markov transform)  which is naturally defined on the
Klein-Dirac quadric. It has been introduced and studied in
\cite{kounchevrenderPade}, \cite{kounchevrenderHIROSHIMA},
\cite{kounchevrenderIsrael}. We will explain briefly this notion here.

First, a crucial point of the pseudo-positive Moment Problem is the remarkable
representation of multivariate polynomials which is sometimes associated with
the names of Gauss and Almansi (cf. \cite{aron}, \cite{sobolev}):
\begin{equation}
P\left(  x\right)  =%
{\displaystyle\sum_{k,\ell}}
p_{k,\ell}\left(  r^{2}\right)  r^{k}Y_{k,\ell}\left(  \theta\right)
,\qquad\text{for }x\in\mathbb{R}^{n},\text{ }r=\left\vert x\right\vert
,\ x=r\theta,\label{Almansi}%
\end{equation}
for $n\in\mathbb{N}$ with $n\geq2.$ Here, for every fixed integer $k\geq0,$
the set
\begin{equation}
\left\{  Y_{k,\ell}\left(  \theta\right)  :\ell=1,2,...,d_{k}\right\}
,\label{Ykl}%
\end{equation}
with
\[
d_{k}=\frac{\left(  2k+n-2\right)  \left(  n+k-3\right)  !}{\left(
n-2\right)  !k!}\qquad\text{for }k\geq0,
\]
is a basis of the space $H_{k}\left(  \mathbb{S}^{n-1}\right)  $ of spherical
harmonics on the unit sphere $\mathbb{S}^{n-1}$ which are homogeneous of
degree $k,$ and $p_{k,\ell}$ are polynomials. It exhibits the representation
of a multivariate polynomial $P$ through the one-dimensional polynomials
$p_{k,\ell}$, having control on the degree: namely, \ if $\Delta^{N}P\left(
x\right)  =0,$ with $\Delta$ being the Laplacian in $\mathbb{R}^{n},$ then
$\deg p_{k,\ell}\leq N-1.$ This representation is alternative to the standard
\[
P\left(  x\right)  =%
{\displaystyle\sum_{\alpha\in\mathbb{Z}_{+}^{n}}}
a_{\alpha}x^{\alpha},
\]
where $\alpha=\left(  \alpha_{1},...,\alpha_{n}\right)  ,$ $\alpha_{i}%
\in\mathbb{N}_{0},$ is a multi-index, $x^{\alpha}=x_{1}^{\alpha_{1}}%
x_{2}^{\alpha_{2}}\cdot\cdot\cdot x_{n}^{\alpha_{n}}.$ It is based on the fact
that the polynomials
\[
r^{k+2j}Y_{k,\ell}\left(  \theta\right)  =\left\vert x\right\vert
^{2j}Y_{k,\ell}\left(  x\right)  \qquad\text{for }j\geq0,\text{ and all
indices }\left(  k,\ell\right)  ,
\]
with $k\geq0$ and $\ell=1,2,...,d_{k},$ are an alternative basis of the space
of all polynomials in $\mathbb{R}^{n}.$ It provides also an alternative
representation for the Moment Problem by considering the moments
\begin{equation}
c_{k,\ell;j}:=%
{\displaystyle\int_{\mathbb{R}^{n}}}
r^{k+2j}Y_{k,\ell}\left(  \theta\right)  d\mu\left(  x\right)  =%
{\displaystyle\int_{\mathbb{R}^{n}}}
\left\vert x\right\vert ^{2j}Y_{k,\ell}\left(  x\right)  d\mu\left(  x\right)
\qquad\text{for }k\geq0,\ \ell=1,2,...,d_{k},\ j\geq
0,\label{MomentProblemMULTI}%
\end{equation}
instead of the usual approach where one considers the Moment Problem in the
form
\[
d_{\alpha}:=%
{\displaystyle\int}
x^{\alpha}d\mu\left(  x\right)  ,\qquad\text{for }\alpha\in\mathbb{N}_{0}.
\]

Further, we consider the \textbf{Hua-Aronszajn} \textbf{kernel} defined by
(cf. the details in \cite{kounchevrenderPade}, \cite{kounchevrenderHIROSHIMA}%
)
\[
K\left(  \zeta\theta;x\right)  :=\frac{\zeta^{n-1}}{r\left(  \zeta
\theta-x\right)  ^{n}},\qquad\zeta\in\mathbb{C},\ \theta\in\mathbb{S}%
^{n-1},\ x\in\mathbb{C}^{n},\ r=\left\vert x\right\vert ,
\]
where $n\in\mathbb{N}$ and $n\geq2.$ Here we use the following notation of
Aronszajn \cite{aron}
\begin{equation}
r\left(  \zeta\theta-x\right)  ^{n}:=\left(  \sqrt{\zeta^{2}-2\zeta
\left\langle \theta,x\right\rangle +\left\vert x\right\vert ^{2}}\right)
^{n},\label{rNotation}%
\end{equation}
where an appropriate root is chosen, $\left\langle \cdot,\cdot\right\rangle $
denoting the scalar product in $\mathbb{C}^{n}$ and $\left\vert \cdot
\right\vert $ the corresponding norm. We notice that the zeros $\zeta
^{2}-2\zeta\left\langle \theta,x\right\rangle +\left\vert x\right\vert ^{2}=0$
are given by
\[
\zeta_{1,2}\left(  \theta,x\right)  =\left\langle \theta,x\right\rangle \pm
i\sqrt{\left\vert x\right\vert ^{2}-\left\vert \left\langle \theta
,x\right\rangle \right\vert ^{2}}.
\]
One has obviously
\begin{equation}
\left\vert \zeta_{1,2}\right\vert =\left\vert x\right\vert .\label{roots}%
\end{equation}
Let us consider the above kernel $K$ for arguments $\left(  \zeta
\theta;x\right)  $ satisfying
\[
\zeta\neq\zeta_{1,2}\left(  \theta,x\right)  .
\]

For us the following representation of the\emph{ Hua-Aronszajn} kernel will be
important (cf. e.g. formula (18) in \cite{kounchevrenderHIROSHIMA}):
\begin{equation}
K\left(  \zeta\theta;x\right)  =\frac{\zeta^{n-1}}{r\left(  \zeta
\theta-x\right)  ^{n}}=\frac{\zeta}{\zeta^{2}-\left\vert x\right\vert ^{2}}%
{\displaystyle\sum_{k,\ell}}
\zeta^{-k}Y_{k,\ell}\left(  \theta\right)  Y_{k,\ell}\left(  x\right)
,\label{Hua-Aronszajn-Expansion}%
\end{equation}
where $\zeta\in\mathbb{C},$ $\theta\in\mathbb{S}^{n-1}$ and $x\in
\mathbb{R}^{n}.$ It is obtained through the representation of $K$ as a product
of the Cauchy kernel and the "complexified Poisson kernel" $\frac{\zeta
^{n-2}\left(  \zeta^{2}-x^{2}\right)  }{r\left(  \zeta\theta-x\right)  ^{n}}$
by the representation:
\[
\frac{\zeta^{n-1}}{r\left(  \zeta\theta-x\right)  ^{n}}=\frac{\zeta}{\zeta
^{2}-x^{2}}\frac{\zeta^{n-2}\left(  \zeta^{2}-x^{2}\right)  }{r\left(
\zeta\theta-x\right)  ^{n}}.
\]
The expansion in spherical harmonics follows directly from the expansion in
spherical harmonics of the "complexified Poisson kernel" (see
\cite{kounchevrenderHIROSHIMA}), \cite{steinWeiss}. 

One may prove now easily the following analog and generalization to the
\textbf{Cauchy formula}:

\begin{proposition}
\label{PAronszajn-Cauchy}For every polynomial $P\left(  x\right)  ,$
$x\in\mathbb{C}^{n},$ the following representation holds:
\begin{equation}
P\left(  x\right)  =\frac{1}{2\pi i}%
{\displaystyle\int_{\mathbb{S}^{1}}}
{\displaystyle\int_{\mathbb{S}^{n-1}}}
\frac{\zeta^{n-1}}{r\left(  \zeta\theta-x\right)  ^{n}}P\left(  \zeta
\theta\right)  d\zeta d\theta\qquad\text{for }x\in\mathbb{C}^{n}.
\label{Aronszajn-Cauchy}%
\end{equation}

\end{proposition}

%

\proof
The proof follows directly from the Almansi representation (\ref{Almansi}),
(\ref{Hua-Aronszajn-Expansion}), and  the classical Cauchy formula, cf.
\cite{Hua}, \cite{aron}.
\endproof

Proposition \ref{PAronszajn-Cauchy} hints at a notion of
\textbf{multidimensional holomorphic functions} on the Klein-Dirac quadric
$\operatorname{KDQ}.$

\begin{definition}
\label{Dholomorphic} We say that the function $f\left(  \zeta,\theta\right)  $
is $\operatorname{KDQ}$-holomorphic on the compact ball $B_{R}\subset
\operatorname{KDQ},$ if it is representable by the series
\begin{equation}
f\left(  \zeta\theta\right)  =%
{\displaystyle\sum_{k,\ell}}
f_{k,\ell}\left(  \zeta^{2}\right)  \zeta^{k}Y_{k,\ell}\left(  \theta\right)
,\qquad\text{for }\zeta\theta\in D, \label{holomorphic}%
\end{equation}
which is absolutely convergent on compacts $K\subset B_{R},$ and where all
functions $f_{k,\ell}$ are analytic in the one-dimensional disc $R\cdot
\mathbb{D},$ where $\mathbb{D}$ denotes the unit disc in $\mathbb{C}.$ In a
similar way, we say that the function $f\left(  \zeta,\theta\right)  $ is
$KDQ-$holomorphic in the \textbf{exterior} of the compact ball $B_{R},$ if $f$
is representable by the series
\begin{equation}
f\left(  \zeta\theta\right)  =%
{\displaystyle\sum_{k,\ell}}
f_{k,\ell}\left(  \zeta^{2}\right)  \frac{1}{\zeta^{k}}Y_{k,\ell}\left(
\theta\right)  ,\qquad\text{for }\zeta\theta\in D, \label{holomorphic2}%
\end{equation}
which is absolutely convergent on compacts $K\subset\operatorname{KDQ}%
\setminus\overline{B_{R}}.$ Here $\overline{B_{R}}$ denotes the topological
closure of $B_{R}.$
\end{definition}

\begin{remark}
The topological boundary $\partial\left(  B_{1}\right)  $ of $B_{1}$ is a
natural boundary for the space of $\operatorname{KDQ}$-holomorphic functions
in $B_{1}.$ In particular, if a function $f$ is $KDQ$-holomorphic and
satisfies
\[%
{\displaystyle\sum_{k,\ell}}
{\displaystyle\int_{0}^{2\pi}}
\left\vert f_{k,\ell}\left(  \zeta^{2}\right)  \right\vert ^{2}d\varphi
<\infty,\qquad\text{for }\ \zeta=e^{i\varphi},
\]
then formula (\ref{Aronszajn-Cauchy}) holds where $P$ is replaced by $f.$ Let
us remark that $\partial\left(  B_{1}\right)  $ is also the Shilov boundary of
the so-called \textbf{fourth domain of Cartan} (or "Lie ball"), and formula
(\ref{Aronszajn-Cauchy}) has been proved in this context for holomorphic
functions of several variables in the Lie ball in $\mathbb{C}^{n},$ (see
\cite{aron}, p. 125, Corollary $1.1$).
\end{remark}

Let $d\mu\left(  x\right)  $ be a signed measure defined on $\mathbb{R}^{n}$
or on some domain $D\subset\mathbb{R}^{n}.$

\begin{definition}
\label{DMarkovStieltjesTransform} We define the (\textbf{multidimensional)
Stieltjes-Markov transform} as the following function on $KDQ,$ whenever the
integral is well defined:
\begin{equation}
\widehat{\mu}\left(  \zeta,\theta\right)  :=%
{\displaystyle\int}
\frac{\zeta^{n-1}d\mu\left(  x\right)  }{r\left(  \zeta\theta-x\right)  ^{n}%
}\qquad\text{for }\ \zeta\theta\in KDQ.\label{Multi-Stieltjes-Markov}%
\end{equation}

\end{definition}

\begin{remark}
In Proposition \ref{PgrowthCondition} below we shall provide a condition under
which the integral on the right hand side of (\ref{Multi-Stieltjes-Markov}) is
absolutely convergent for all $\operatorname{Im}\zeta^{2}>0,$ for $\left\vert
\zeta\right\vert $ sufficiently large.
\end{remark}

Further, let us introduce the following important notion of \textbf{component
measures, }see \cite{kounchevrenderARXIV} or \cite{kounchevrender}:

\begin{definition}
\label{Dcomponent} For every index $\left(  k,\ell\right)  ,$ $k\in
\mathbb{N}_{0},$ $\ell=1,2,...,d_{k},$ the one-dimensional \textbf{component
measure} $\mu_{k,\ell}\left(  r\right)  $ on $\mathbb{R}_{+}$ is uniquely
defined as the measure on $\mathbb{R}_{+},$ such that for every continuous
function $F$ with a compact support in $\mathbb{R}_{+}$ it  satisfies the
following identity:
\[%
{\displaystyle\int_{0}^{\infty}}
F\left(  r\right)  d\mu_{k,\ell}\left(  r\right)  =%
{\displaystyle\int_{\mathbb{R}^{n}}}
F\left(  r\right)  Y_{k,\ell}\left(  \theta\right)  d\mu\left(  x\right)  ,
\]
where $x=r\theta,$ $r=\left\vert x\right\vert ,$ $\theta\in\mathbb{S}^{n-1}.$
Formally, we will write
\[
d\mu_{k,\ell}\left(  r\right)  :=\int_{\mathbb{S}_{\theta}^{n-1}}Y_{k,\ell
}\left(  \theta\right)  d\mu\left(  r\theta\right)  .
\]

\end{definition}

The transform $\widehat{\mu}\left(  \zeta,\theta\right)  $ will be interesting
for us only for large values of $\left\vert \zeta\right\vert .$ However, let
us note that a priori it is not clear for which values of the parameters
$\zeta,$ $\theta$ does the function $\widehat{\mu}\left(  \zeta,\theta\right)
$ make sense, i.e. the integral in (\ref{Multi-Stieltjes-Markov}) is well
defined. If we assume for simplicity that the signed measure $\mu$ has a
compact support in the ball $B_{R}$ in $\mathbb{R}^{n}$ for some $R>0,$ then
by (\ref{roots}) $\widehat{\mu}\left(  \zeta,\theta\right)  $ is well defined
for all $\zeta$ with $\left\vert \zeta\right\vert >R.$ But for a general
measure singularities can arise for all possible values of $\zeta$ as far as
$\left\vert \zeta\right\vert =\left\vert x\right\vert $ on the support of
$\mu.$ The following representation is \emph{remarkable} since it shows that
the multidimensional Stieltjes-Markov transform $\widehat{\mu}\left(
\zeta,\theta\right)  $ is representable through infinitely many
one-dimensional Stieltjes transforms, those of the measures $r^{k}d\mu
_{k,\ell}\left(  r\right)  .$ More precisely, we have:

\begin{proposition}
\label{PSeries-Stieltjes-Transform} Let the signed measure $\mu$ have a
compact support in the ball $B_{R}\subset\mathbb{R}^{n}$ and have a finite
variation there. Then the following expansion
\begin{equation}
\widehat{\mu}\left(  \zeta,\theta\right)  =%
{\displaystyle\sum_{k,\ell}}
\zeta^{-k+1}Y_{k,\ell}\left(  \theta\right)  \times%
{\displaystyle\int_{0}^{R}}
\frac{1}{\zeta^{2}-r^{2}}r^{k}d\mu_{k,\ell}\left(  r\right)  \qquad\text{for
all }\left\vert \zeta\right\vert >R,\ \theta\in\mathbb{S}^{n-1}
\label{Series-Stieltjes-Markov}%
\end{equation}
is absolutely convergent on compacts in $\operatorname{KDQ}\setminus
\overline{B_{R}}.$
\end{proposition}

%

\proof
First, we have the representation formula (\ref{Hua-Aronszajn-Expansion}), and
we use the definition of the component measures $\mu_{k,\ell},$ cf. Definition
\ref{Dcomponent}, to obtain the following:\
\begin{align*}
\widehat{\mu}\left(  \zeta,\theta\right)   &  =%
{\displaystyle\sum_{k,\ell}}
\zeta^{-k}Y_{k,\ell}\left(  \theta\right)
{\displaystyle\int_{\mathbb{R}^{n}}}
Y_{k,\ell}\left(  x\right)  \frac{\zeta}{\zeta^{2}-\left\vert x\right\vert
^{2}}d\mu\left(  x\right) \\
&  =%
{\displaystyle\sum_{k,\ell}}
\zeta^{-k+1}Y_{k,\ell}\left(  \theta\right)
{\displaystyle\int_{0}^{R}}
\frac{\zeta}{\zeta^{2}-r^{2}}r^{k}d\mu_{k,\ell}\left(  r\right)  .
\end{align*}
Its absolute convergence on compacts $K\subset\operatorname{KDQ}%
\setminus\overline{B_{R}}$ is obtained by a routine estimation using the bound
$k^{n-2}$ for the spherical harmonics $Y_{k,\ell}\left(  \theta\right)  ,$ cf.
\cite{steinWeiss}.%

\endproof

At this point we define the very important class of pseudo-positive measures:

\begin{definition}
We will say that the measure $\mu$ is \textbf{pseudo-positive} in the ball
$B_{R}$ in $\operatorname{KDQ}$ (with respect to the fixed basis $\left\{
Y_{k,\ell}\left(  \theta\right)  :k\in\mathbb{N}_{0},\ \ell=1,2,...,d_{k}%
\right\}  $ ) if it satisfies
\[
d\mu_{k,\ell}\left(  r\right)  \geq0\qquad\text{for  }0\leq r\leq R,
\]
for all indices $\left(  k,\ell\right)  ,$ (see Definition \ref{Dcomponent}).
\end{definition}

We have the following simple condition for convergence of the multidimensional
Stieltjes-Markov transform:

\begin{proposition}
\label{PgrowthCondition} Assume that for the signed measure $\mu,$ and for
some constants $C>0$ and $D>0,$ its components $\mu_{k,\ell}$ satisfy the
following \textbf{growth condition}%
\[
\left\vert
{\displaystyle\int_{0}^{\infty}}
r^{k}d\mu_{k,\ell}\left(  r\right)  \right\vert \leq CD^{k}\qquad\text{for all
}\left(  k,\ell\right)  \text{ with }k\geq0,
\]
for some constants $C,D>0,$ independent of $k,\ell.$ Then it follows that the
multidimensional Stieltjes-Markov transform $\widehat{\mu}\left(  \zeta
,\theta\right)  $ in  Definition \ref{DMarkovStieltjesTransform} is absolutely
convergent for all $\zeta$ with $\operatorname{Im}\zeta^{2}>0,$ $\left\vert
\zeta\right\vert >D$ and all $\theta\in\mathbb{S}^{n-1}.$
\end{proposition}

%

\proof
The proof follows by a direct estimate of the series
(\ref{Series-Stieltjes-Markov}) taking into account the bound $k^{n-2}$ for
the spherical harmonics $Y_{k,\ell}\left(  \theta\right)  ,$ cf.
\cite{steinWeiss}.%

\endproof

Proposition \ref{PgrowthCondition} is remarkable in some sense, since it shows
that for a reasonable class of measures with \textbf{non-compact support}, one
has convergence of the Stieltjes-Markov transform $\widehat{\mu}$ for all
sufficiently large values of $\zeta,$ thanks to the expansion
(\ref{Series-Stieltjes-Markov}).

\section{Nevanlinna class on the Klein-Dirac quadric and generalized
Hamburger-Nevanlinna theorem \label{SNevan}}

In the present section we will discuss the role of the multidimensional
Stieltjes-Markov transform for the multidimensional pseudo-positive Moment Problem.

Let us recall the classical Nevanlinna class of analytical functions in the
half-plane: A function $f$ defined in the upper half plane $\mathbb{C}$
belongs to the Nevanlinna class $\mathcal{N}$ iff it satisfies
\[
\operatorname{Im}z>0\text{ implies }\operatorname{Im}f\left(  z\right)
\geq0.
\]
The elements $f$ of $\mathcal{N}$ have the following
\textbf{Riesz-Herglotz-Nevanlinna} representation (\cite{Akhiezer}, p. $92$):
\[
f\left(  z\right)  =az+b+%
{\displaystyle\int_{-\infty}^{\infty}}
\frac{1+uz}{u-z}d\nu\left(  u\right)  ,\qquad\text{for }z\in\mathbb{C}%
,\ \operatorname{Im}z\neq0,
\]
for some non-negative measure $\nu$ with bounded variation (equal to
$\operatorname{Im}f\left(  i\right)  -a$), and real constants $a\geq0$ and
$b.$

We recall the classical \textbf{Hamburger-Nevanlinna} theorem which plays an
important role in the one-dimensional Moment Problem (cf. Theorem $3.2.1$ in
\cite{Akhiezer}). It shows that the Stieltjes transform of a measure encodes
in a very appropriate way the information about the moments of the measure,
and the moments may be extracted only by means of an adequate "asymptotic
expansion" and not by the usual Taylor expansion, which is asymptotic but in
general notoriously divergent.

\begin{theorem}
\label{THamburger-Nevanlinna} If the non-negative measure $\sigma$ solves the
truncated Moment Problem
\begin{equation}
s_{j}=%
{\displaystyle\int_{-\infty}^{\infty}}
u^{j}d\sigma\left(  u\right)  \qquad\text{for }%
j=0,1,...,2N,\label{Sigma-moments}%
\end{equation}
then the Stieltjes transform of the measure $\sigma$ given by
\begin{equation}
f\left(  z\right)  =%
{\displaystyle\int_{-\infty}^{\infty}}
\frac{d\sigma\left(  u\right)  }{u-z}\label{Stieltjes2}%
\end{equation}
belongs to the Nevanlinna class $\mathcal{N}$ and for every fixed $\delta>0$ (
$<\pi/2$ ) the following limit holds uniformly for $\delta\leq\arg z\leq
\pi-\delta$ :
\begin{equation}
\lim_{z\longrightarrow\infty}z^{2N+1}\left\{  f\left(  z\right)  +%
{\displaystyle\sum_{j=0}^{2N-1}}
\frac{s_{j}}{z^{j+1}}\right\}  =-s_{2N}.\label{Nevanlinna-limit}%
\end{equation}
In particular, the formal series $%
{\displaystyle\sum_{j=0}^{\infty}}
\frac{s_{j}}{z^{j+1}}$ gives an asymptotic expansion to all orders for
$f\left(  z\right)  $ as $z\longrightarrow\infty.$

Conversely, if for some function $f\in\mathcal{N}$ equation
(\ref{Nevanlinna-limit}) holds with $s_{j}\in\mathbb{R},$ for all
$N\in\mathbb{N}$, at least for $z=iy,$ $y\longrightarrow\infty,$ then $f$ is
representable in the form (\ref{Stieltjes2}) and the non-negative measure
$\sigma$ has finite moments of any order, and in particular satisfies
(\ref{Sigma-moments}).
\end{theorem}

We define the multidimensional pseudo-positive Nevanlinna class on the
Klein-Dirac quadric:

\begin{definition}
Let the function $f$ be $\operatorname{KDQ}$-holomorphic in the exterior of
the ball $B_{R}\subset\operatorname{KDQ},$ in the sense of Definition
\ref{Dholomorphic}; hence $f$ is representable by the series
(\ref{holomorphic2}). We say that $f$ belongs to the \textbf{pseudo-positive
Nevanlinna} class on $\operatorname{KDQ}\setminus B_{R}$ iff all its
components $f_{k,\ell}$ belong to the one-dimensional Nevanlinna class
$\mathcal{N}.$
\end{definition}

Now we are able to prove the generalization of the Hamburger-Nevanlinna
theorem which shows that the multidimensional Stieltjes-Markov transform
encodes in an asymptotic expansion form the information about the
multidimensional moments of a measure (cf. \cite{kounchevrender}). However,
since the component measures $\mu_{k,\ell}$ are not so strongly correlated we
will need to keep them tight by means of some growth condition.

\begin{theorem}
\label{THamburger-Nevanlinna-MULTI} Assume that a \textbf{pseudo-positive
measure} $\mu$ is given, and solves the following multidimensional Moment
Problem, (\ref{MomentProblemMULTI}),
\[
s_{k,\ell;j}=%
{\displaystyle\int_{\mathbb{R}^{n}}}
r^{k+2j}Y_{k,\ell}\left(  \theta\right)  d\mu\left(  x\right)  \qquad\text{
for all }\left(  k,\ell\right)  \text{ and }j=0,1,...,2N.
\]
Assume that, for some constants $C>0$ and $D>0,$ independent of $k,$ $\ell,$
the following growth condition holds:
\begin{equation}
\left\vert s_{k,\ell;0}\right\vert \leq CD^{k}\qquad\text{for all indices
}\left(  k,\ell\right)  .\label{growth-condition}%
\end{equation}
Then the multidimensional Stieltjes-Markov transform $f\left(  \zeta
,\theta\right)  =\widehat{\mu}\left(  \zeta,\theta\right)  $ defined by
(\ref{Multi-Stieltjes-Markov}) is absolutely convergent and belongs to the
pseudo-positive Nevanlinna class on $\operatorname{KDQ}\setminus
\overline{B_{D}}.$ For every index $\left(  k,\ell\right)  $ the following
asymptotic expansion holds uniformly for $\delta\leq\arg\zeta^{2}\leq
\pi-\delta,$ and $\zeta\longrightarrow\infty$ (but \textbf{not uniformly} on
the indices $k,\ell$ ! ):
\begin{equation}
\zeta^{4N+1}\left\{  \zeta^{k-1}%
{\displaystyle\int_{\mathbb{S}^{n-1}}}
f\left(  \zeta,\theta\right)  Y_{k,\ell}\left(  \theta\right)  d\theta-%
{\displaystyle\sum_{j=0}^{2N-1}}
\frac{s_{k,\ell;j}}{\zeta^{k+2j}}\right\}  =s_{k,\ell;2N}%
.\label{Asymptotic-Nevanlinna}%
\end{equation}
In particular, the formal series $%
{\displaystyle\sum_{j=0}^{2N-1}}
\frac{s_{k,\ell;j}}{\zeta^{k+2j}}$ provides an asymptotic expansion to all
orders for the component $%
{\displaystyle\int_{\mathbb{S}^{n-1}}}
f\left(  \zeta,\theta\right)  Y_{k,\ell}\left(  \theta\right)  d\theta$ of $f$
along $Y_{k,\ell}\left(  \theta\right)  .$
\end{theorem}

%

\proof
By formula (\ref{Series-Stieltjes-Markov}) we see that
\[
\zeta^{k-1}%
{\displaystyle\int_{\mathbb{S}^{n-1}}}
\widehat{\mu}\left(  \zeta,\theta\right)  Y_{k,\ell}\left(  \theta\right)
d\theta=%
{\displaystyle\int_{0}^{\infty}}
\frac{1}{\zeta^{2}-r^{2}}r^{k}d\mu_{k,\ell}\left(  r\right)  ,
\]
and the right hand side is a one-dimensional Stieltjes transforms. It follows
by Theorem \ref{THamburger-Nevanlinna} that for every $\left(  k,\ell\right)
$ the component measures $\mu_{k,\ell}$ (which are non-negative measures) have
the moments
\[
s_{k,\ell;j}=%
{\displaystyle\int_{0}^{\infty}}
r^{k+2j}d\mu_{k,\ell}\left(  r\right)  \qquad\text{for }j=0,1,...,2N,
\]
and the asymptotic expansion (\ref{Asymptotic-Nevanlinna}) holds as in the
classical case covered by Theorem \ref{THamburger-Nevanlinna}.%

\endproof

There is a converse statement of the above Theorem
\ref{THamburger-Nevanlinna-MULTI} which we do not formulate.

\begin{remark}
In the one-dimensional case the Moment Problem provides one of the most
important methods for summation of divergent series, cf. \cite{hardy},
\cite{jones-thron}, chapter $9.$ Hence, Theorem
\ref{THamburger-Nevanlinna-MULTI} generalizes known classical results about
summation of divergent series. Indeed, in the classical case, by Theorem
\ref{THamburger-Nevanlinna} we attach the function $f\left(  z\right)  $ to
the divergent series $\sum\frac{s_{j}}{z^{j+1}}$ using the fact that the
sequence $\left\{  s_{j}\right\}  $ is positive-definite. Now, in a similar
manner,  Theorem \ref{THamburger-Nevanlinna-MULTI} provides a new summation
method for multidimensional \textbf{divergent series}: The multidimensional
Stieltjes-Markov transform $\widehat{\mu}\left(  \zeta,\theta\right)  $ is an
asymptotic expansion in polyharmonic functions, associated with the formal
Laurent series (divergent in general!)
\[
f_{N}\left(  \zeta\theta\right)  =\frac{1}{\zeta}%
{\displaystyle\sum_{k,\ell}}
{\displaystyle\sum_{j=0}^{2N-1}}
\frac{s_{k,\ell;j}}{\zeta^{k+2j}}Y_{k,\ell}\left(  \theta\right)  ,
\]
which satisfies formally the polyharmonic equation  $\Delta^{2N}f\left(
\zeta\theta\right)  =0,$ i.e. we may write formally%
\[
\lim_{\zeta\longrightarrow\infty}\zeta^{4N+1}\left\{  \widehat{\mu}\left(
\zeta,\theta\right)  -f_{N}\left(  \zeta\theta\right)  \right\}  -g_{N}\left(
\zeta,\theta\right)  =0
\]
where we have put, again formally,
\[
g_{N}\left(  \zeta,\theta\right)  =%
{\displaystyle\sum_{k,\ell}}
\frac{Y_{k,\ell}\left(  \theta\right)  }{\zeta^{k}}c_{k,\ell;2N};
\]
note that the function $g_{N}$ satisfies formally the polyharmonic equation
$\Delta^{2N+1}g_{N}=0.$
\end{remark}

\section{Multidimensional Isospectral Deformation and pseudo-positive Toda
lattice\label{Smulti}}

We are about to define our multidimensional generalization of Toda lattice. We
proceed by "reverse engineering": we first define the Toda lattice in the
model variables $\lambda$ and $r$  of the multidimensional Moment Problem,
and then by reversion we come to a formulation in the "physical" variables $x$
and $y.$ 

We will consider the measures $d\mu\left(  x\right)  $ for $x=r\theta
\in\mathbb{R}^{n}$ and $r=\left\vert x\right\vert $ which have the following
Laplace-Fourier expansion%
\begin{equation}
d\mu\left(  x\right)  =%
{\displaystyle\sum_{k,\ell}}
d\mu_{k,\ell}\left(  r\right)  Y_{k,\ell}\left(  \theta\right)  \qquad
\text{for }k\in\mathbb{N}_{0},\ \ell=1,...,d_{k}. \label{dmuNew}%
\end{equation}

We have seen that in the one-dimensional case the \textbf{finite Toda lattice}
is generated by an "isospectral deformation" of a non-negative measure
$\mu\left(  t\right)  $ and the identity (\ref{CF}), namely by
\[
f\left(  \lambda,t\right)  =%
{\displaystyle\int}
\frac{d\mu\left(  \tau;t\right)  }{\lambda-\tau}=\frac{1}{\lambda-b_{N}\left(
t\right)  -\frac{a_{N-1}^{2}\left(  t\right)  }{\lambda-b_{N-1}-\cdot
\cdot\cdot}}=\frac{Q_{N}\left(  \lambda,t\right)  }{P_{N}\left(
\lambda,t\right)  }=\sum_{j=1}^{N}\frac{r_{j}^{2}\left(  t\right)  }%
{\lambda-\lambda_{j}},
\]
for $\lambda\notin\sigma\left(  L\right)  ,$ $t\in\mathbb{R}.$ The spectrum of
the measure $d\mu$ coincides with the set of the constants $\lambda_{j},$ cf.
\cite{nikishinSorokin}. "Isospectrality" is understood as time-independence of
the spectrum of the measure $d\mu\left(  t\right)  .$ 

In a similar way, we propose to generalize the "isospectral property" by
considering a time-dependent pseudo-positive measure $\mu$ and its
multidimensional Stieltjes-Markov transform $\widehat{\mu}\left(  \zeta
,\theta;t\right)  $ on the Klein-Dirac quadric.  We will have respectively the
series
\begin{align}
f\left(  \zeta,\theta;t\right)   &  =\widehat{\mu}\left(  \zeta,\theta
;t\right)  =%
{\displaystyle\int}
\frac{\zeta^{n-1}d\mu\left(  x;t\right)  }{r\left(  \zeta\theta-x\right)
^{n}}=%
{\displaystyle\sum_{k,\ell}}
\frac{Y_{k,\ell}\left(  \theta\right)  }{\zeta^{k}}%
{\displaystyle\int_{0}^{\infty}}
\frac{\zeta}{\zeta^{2}-r^{2}}r^{k}d\mu_{k,\ell}\left(  r;t\right)
\label{fztheta}\\
&  =%
{\displaystyle\sum_{k,\ell}}
\frac{Y_{k,\ell}\left(  \theta\right)  }{\zeta^{k-1}}%
{\displaystyle\sum_{j=1}^{N}}
\frac{\lambda_{k,\ell;j}^{k}r_{k,\ell;j}^{2}\left(  t\right)  }{\zeta
^{2}-\lambda_{k,\ell;j}^{2}}\qquad\qquad\zeta\in\mathbb{C},\ \theta
\in\mathbb{S}^{n-1},\nonumber
\end{align}
which converges in the absolute sense. We define the spectrum of the measure
$d\mu$ as the union of the spectral of all measures $\left\{  \mu_{k,\ell
}\right\}  _{k,\ell}$ which is the set of all constants $\left\{
\lambda_{k,\ell;j}\right\}  .$ Now, by analogy with the one-dimensional case,
under  "isospectrality" we will understand the time-independence of all
constants $\lambda_{k,\ell;j}.$

We will consider the more convenient measures $d\widetilde{\mu}_{k,\ell
}\left(  \rho;t\right)  $ where we have put $\rho=r^{2},$ hence,
\begin{equation}
d\widetilde{\mu}_{k,\ell}\left(  \rho;t\right)  =r^{k}d\mu_{k,\ell}\left(
r;t\right)  =\rho^{k/2}d\mu_{k,\ell}\left(  \sqrt{\rho};t\right)
,\label{dmuTilde}%
\end{equation}
i.e.
\[
d\widetilde{\mu}_{k,\ell}\left(  \rho;t\right)  =%
{\displaystyle\sum_{j=1}^{N}}
r_{k,\ell;j}^{2}\left(  t\right)  \lambda_{k,\ell;j}^{k}\delta\left(
\rho-\lambda_{k,\ell;j}^{2}\right)  .
\]
We put
\begin{align}
\widetilde{r}_{k,\ell;j}^{2}\left(  t\right)   &  :=r_{k,\ell;j}^{2}\left(
t\right)  \lambda_{k,\ell;j}^{k}\label{rtilde}\\
\widetilde{\lambda}_{k,\ell;j} &  :=\lambda_{k,\ell;j}^{2}.\nonumber
\end{align}
By the definition of the measure $\mu_{k,\ell}$ we obtain the equality
\[
f\left(  \zeta,\theta;t\right)  =%
{\displaystyle\sum_{k,\ell}}
\frac{Y_{k,\ell}\left(  \theta\right)  }{\zeta^{k-1}}%
{\displaystyle\int_{0}^{\infty}}
\frac{1}{\zeta^{2}-\rho}d\widetilde{\mu}_{k,\ell}\left(  \rho;t\right)
\]
and the Jacobi matrix for $d\widetilde{\mu}_{k,\ell}\left(  \rho;t\right)  $
will satisfy a one-dimensional Toda lattice equation, called below
\emph{pseudo-positive Toda lattice}. We now give the rigorous definition.

\begin{definition}
\label{DpseudoToda} We say that the pseudo-positive measure $\mu\left(
x;t\right)  ,$ $x\in\mathbb{R}^{n},$ given for every $t\in\mathbb{R}$ by
(\ref{dmuNew}), with $\mu_{k,\ell}\left(  r\right)  $ replaced by $\mu
_{k,\ell}\left(  r;t\right)  ,$ represents a \textbf{pseudo-positive Toda
lattice} if all component measures $\mu_{k,\ell}$ are given by
\[
d\mu_{k,\ell}\left(  r;t\right)  =%
{\displaystyle\sum_{j}}
r_{k,\ell;j}^{2}\left(  t\right)  \delta\left(  r-\lambda_{k,\ell;j}\right)
dr,\qquad r\in\left[  0,\infty\right)  ,
\]
where for all indices we have
\[
r_{k,\ell;j}\geq0,\quad\lambda_{k,\ell;j}\geq0.
\]
We assume further that the system is isospectral, i.e. the spectrum $\left\{
\lambda_{k,\ell;j}\right\}  _{j=1}^{N}$ of all measures $\mu_{k,\ell}\left(
r;t\right)  $ is constant, i.e. independent of $t\in\mathbb{R}.$ We assume
that the masses $\widetilde{r}_{k,\ell;j}\left(  t\right)  $ given by
(\ref{rtilde}) vary with the time $t$ according to the Toda lattice equations
(\ref{TodaSolution})-(\ref{TodaSolution2}), i.e.
\begin{align}
\widetilde{\lambda}_{k,\ell;j}\left(  t\right)   &  =\widetilde{\lambda
}_{k,\ell;j}\left(  0\right)  \label{TodaSolution-kl}\\
\widetilde{r}_{k,\ell;j}^{2}\left(  t\right)   &  =\frac{\widetilde{r}%
_{k,\ell;j}^{2}\left(  0\right)  e^{-2\widetilde{\lambda}_{k,\ell;j}t}}{%
{\displaystyle\sum_{m=1}^{N}}
\widetilde{r}_{k,\ell;m}^{2}\left(  0\right)  e^{-2\widetilde{\lambda}%
_{k,\ell;m}t}},\nonumber
\end{align}
with the restriction
\begin{equation}%
{\displaystyle\sum_{j=1}^{N}}
\widetilde{r}_{k,\ell;j}^{2}\left(  0\right)  =1\qquad\text{for all }\left(
k,\ell\right)  .\label{sumrtilde=1}%
\end{equation}
(see \cite{MoserNotes}, p. $223,$ \cite{simonSzego}).
\end{definition}

From (\ref{TodaSolution-kl}) we obtain immediately the equalities
\begin{align*}
\widetilde{r}_{k,\ell;j}^{2}\left(  t\right)   &  =\frac{\lambda_{k,\ell
;j}^{k}r_{k,\ell;j}^{2}\left(  0\right)  e^{-2\lambda_{k,\ell;j}^{2}t}}{%
{\displaystyle\sum_{m=1}^{N}}
\lambda_{k,\ell;m}^{k}r_{k,\ell;m}^{2}\left(  0\right)  e^{-2\lambda
_{k,\ell;m}^{2}t}}\\
r_{k,\ell;j}^{2}\left(  t\right)   &  =\frac{r_{k,\ell;j}^{2}\left(  0\right)
e^{-2\lambda_{k,\ell;j}^{2}t}}{%
{\displaystyle\sum_{m=1}^{N}}
\lambda_{k,\ell;m}^{k}r_{k,\ell;m}^{2}\left(  0\right)  e^{-2\lambda
_{k,\ell;m}^{2}t}}%
\end{align*}
and from (\ref{sumrtilde=1})
\begin{equation}%
{\displaystyle\sum_{m=1}^{N}}
\lambda_{k,\ell;m}^{k}r_{k,\ell;m}^{2}\left(  0\right)  =1. \label{sumr=1}%
\end{equation}

Let us define the \textbf{Jacobi matrix} $L_{k,\ell}$ related  to the
Stieltjes transform%
\[%
{\displaystyle\int_{0}^{\infty}}
\frac{1}{z-\rho}d\widetilde{\mu}_{k,\ell}\left(  \rho\right)  \qquad\text{for
}\operatorname{Im}z\neq0,
\]
by putting
\[
L_{k,\ell}\left(  t\right)  :=\left(  \left(  \widetilde{a}_{k,\ell;j}\left(
t\right)  ,\ \widetilde{b}_{k,\ell;j}\left(  t\right)  \right)  \right)
_{j=0}^{N};
\]
it generalizes  the one-dimensional Jacobi matrix $L$ in (\ref{L}), cf.
\cite{moser1}, p. $473,$ also \cite{simon98}. Now we are able to define the
Hamiltonians corresponding to the system of equations (\ref{traceL2}), by
putting
\begin{equation}
\frac{1}{2}H_{k,\ell}=2\left\{
{\displaystyle\sum_{j=1}^{N-1}}
\widetilde{a}_{k,\ell;j}^{2}+\frac{1}{2}%
{\displaystyle\sum_{j=1}^{N}}
\widetilde{b}_{k,\ell;j}^{2}\right\}  =%
{\displaystyle\sum_{j=1}^{N}}
\widetilde{\lambda}_{k,\ell;j}^{2}=%
{\displaystyle\sum_{j=1}^{N}}
\lambda_{k,\ell;j}^{4}\label{Htilde}%
\end{equation}
The corresponding total Hamiltonian is now given by
\begin{equation}
H=%
{\displaystyle\sum_{k,\ell}}
H_{k,\ell}=2%
{\displaystyle\sum_{k,\ell}}
{\displaystyle\sum_{j=1}^{N}}
\widetilde{\lambda}_{k,\ell;j}^{2}=2%
{\displaystyle\sum_{k,\ell}}
{\displaystyle\sum_{j=1}^{N}}
\lambda_{k,\ell;j}^{4}\label{HamiltonMulti}%
\end{equation}

We obtain 

\begin{corollary}
\label{CorFlaschkaMULTI}1. The variables $\left(  \widetilde{a}_{k,\ell
;j},\widetilde{b}_{k,\ell;j}\right)  _{k,\ell;j}$ satisfy equations
\begin{align}
a_{k,\ell;j}^{\prime} &  =a_{k,\ell;j}\left(  b_{k,\ell;j+1}-b_{k,\ell
;j}\right)  ,\qquad j=1,2,...,N-1\label{ab1MULTI}\\
b_{k,\ell;j}^{\prime} &  =2\left(  a_{k,\ell;j}^{2}-a_{k,\ell;j-1}^{2}\right)
,\qquad\text{ }j=1,2,...,N\label{ab2MULTI}\\
a_{k,\ell;0}^{2} &  =0,\quad a_{k,\ell;N}^{2}=0.\label{ab3MULTI}%
\end{align}
which coincide with equations (\ref{ab1})-(\ref{ab3}). 

2. By using the one-to-one map defined by the equations (\ref{FlaschkaChange}%
)-(\ref{FlaschkaChange3}) we find the variables $\left(  x_{k,\ell
;j},y_{k,\ell;j}\right)  _{k,\ell;j}$ which satisfy the system of equations
(\ref{TodaLattice}). This represents a Hamiltonian system with Hamiltonian
given by (\ref{HamiltonMulti}).
\end{corollary}

From above the following important statement is easily derived, by using
(\ref{rtilde}):

\begin{proposition}
\label{PHamiltonConverge} Let us assume that
\begin{equation}
C_{\lambda}:=%
{\displaystyle\sum_{k,\ell}}
{\displaystyle\sum_{j=1}^{N}}
\lambda_{k,\ell;j}^{4}<\infty.\label{SumLambdas}%
\end{equation}
Then for every $t\geq0$ the Hamiltonian $H$ of the pseudo-positive Toda
lattice is a real number independent of $t\in\mathbb{R}$ which is bounded by
$2C_{\lambda}.$
\end{proposition}

The following Theorem shows that a very mild restriction on the spectrum
$\left\{  \lambda_{k,\ell;j}\right\}  $ suffices to guarantee convergence of
the multidimensional Stieltjes-Markov transform for every $t\geq0.$

\begin{theorem}
\label{Tnice} Let a pseudo-positive Toda lattice be given in the sense of
Definition \ref{DpseudoToda}. Then the multidimensional Stieltjes-Markov
transform $\widehat{\mu}\left(  \zeta,\theta;t\right)  $ defined in
(\ref{fztheta}) (of the associated pseudo-positive measure $\mu\left(
x;t\right)  $ defined by (\ref{dmuNew}) with $\mu_{k,\ell}\left(  r\right)  $
replaced by the measures $\mu_{k,\ell}\left(  r;t\right)  $ in Definition
\ref{DpseudoToda}) is convergent (in the sense, the integral representing it
is convergent) for all $t\geq0$ and for all $\left\vert \zeta\right\vert >1$
with $\operatorname{Im}\zeta^{2}>0.$
\end{theorem}

%

\proof
From (\ref{TodaSolution-kl}) it follows
\begin{align*}%
{\displaystyle\int_{0}^{\infty}}
r^{k}d\mu_{k,\ell}\left(  r\right)   &  =%
{\displaystyle\sum_{j=1}^{N}}
\lambda_{k,\ell;j}^{k}r_{k,\ell;j}^{2}\left(  t\right)  =%
{\displaystyle\sum_{j=1}^{N}}
\lambda_{k,\ell;j}^{k}\frac{r_{k,\ell;j}^{2}\left(  0\right)  e^{-2\lambda
_{k,\ell;j}t}}{%
{\displaystyle\sum_{m=1}^{N}}
\lambda_{k,\ell;m}^{k}r_{k,\ell;m}^{2}\left(  0\right)  e^{-2\lambda
_{k,\ell;m}t}}\\
&  =1.
\end{align*}
The proof is finished by Proposition \ref{PgrowthCondition} since condition
(\ref{growth-condition}) is fulfilled with $C=D=1.$%

\endproof

\section{The functional model for the pseudo-positive Jacobi (Schr\"{o}dinger)
operator \label{SJacobi}}

We will define an operator in a functional model generated by a given
pseudo-positive measure $\mu.$ This model is related to the above
multidimensional Toda lattice in a way generalizing the one-dimensional
considered by J. Moser in \cite{moser1}, \cite{moser2}.

We recall that in the classical spectral theory, for a non-negative measure
$\mu$ on $\mathbb{R},$ the functional model is defined by the space
\[
L_{2}\left(  \mu\right)  =\left\{  f:%
{\displaystyle\int_{\mathbb{R}}}
\left\vert f\left(  t\right)  \right\vert ^{2}d\mu\left(  t\right)
<\infty\right\}
\]
and the multiplication operator $A_{t}f=tf\left(  t\right)  ,$ for
$t\in\mathbb{R}$, cf. \cite{reedSimon}. We assume that the polynomials are
dense in the space $L_{2}\left(  \mu\right)  ,$ i.e. we assume that we have a
determinate Moment Problem relative to the measure $\mu$ (cf. \cite{Akhiezer}%
). The orthonormal basis is given by polynomials $\left\{  P_{j}\left(
t\right)  \right\}  _{j\geq0}$ which are orthonormal with respect to the
measure $\mu.$ In this basis the operator $A_{t}$ is represented by the Jacobi
matrix obtained through the corresponding $3-$term recurrence relations. See
\cite{Akhiezer}, \cite{simon98}, \cite{teschl} for the details.

In the present situation the functional model is constructed as follows: We
assume that a pseudo-positive measure $\mu$ is given. We define the space
\begin{align*}
&
{\textstyle\bigoplus_{k,\ell}^{\prime}}
L_{2}\left(  \mu_{k,\ell}\right)  \\
&  =\left\{  f\equiv\left(  f_{k,\ell}\right)  ,k\in\mathbb{N}_{0}%
,\ell=1,2,...,d_{k}:%
{\displaystyle\sum_{k,\ell}}
{\displaystyle\int_{0}^{\infty}}
\left\vert f_{k,\ell}\left(  r\right)  \right\vert ^{2}d\mu_{k,\ell}\left(
r\right)  <\infty\right\}  ,
\end{align*}
where the measures $\mu_{k,\ell}$ are non-negative component measures on
$\mathbb{R}_{+}$ associated with $\mu$ according to Definition
\ref{Dcomponent}. The operator $A_{\left\vert x\right\vert ^{2}}$ is defined
as the multiplication operator
\[
A_{\left\vert x\right\vert ^{2}}f\left(  x\right)  =\left\vert x\right\vert
^{2}f\left(  x\right)  \qquad\text{for }x\in\mathbb{R}^{n},
\]
where $f\in D\left(  \left\vert \cdot\right\vert ^{2}\right)  \subset
L^{2}\left(  \mathbb{R}^{n}\right)  ;$ here $D$ is the appropriate domain of
definition of the operator $A_{\left\vert x\right\vert ^{2}}.$ Further, we
assume that all Moment Problems defined by the measures $\mu_{k,\ell}$ are
determined on $\mathbb{R}_{+},$ cf. \cite{kounchevrender}. In this case the
basis for the space $%
{\textstyle\bigoplus_{k,\ell}^{\prime}}
L_{2}\left(  \mu_{k,\ell}\right)  $ is the set of polynomials
\[
r^{k}P_{k,\ell;j}\left(  r^{2}\right)
\]
where the polynomials $\left\{  P_{k,\ell;j}\right\}  _{j\geq0}$ are
orthonormal relative to the measure $\widetilde{\mu}_{k,\ell}$ on
$\mathbb{R}_{+},$ defined in terms of the $d\mu_{k,\ell}$ by (\ref{dmuTilde}).
In this basis the operator $A_{\left\vert x\right\vert ^{2}}$ is represented
by a matrix with components given by the Jacobi matrices $L_{k,\ell}$ in the
spaces $L_{2}\left(  \widetilde{\mu}_{k,\ell}\right)  .$ The representation of
the Hamiltonian in (\ref{Htilde}) means that this operator is bounded in the
respective $L_{2}$ space.

\subsection{Physical meaning of the pseudo-positive Toda lattice}

We consider here the pseudo-positive Toda lattice in a generalization of the
Flaschka variables, and we also provide some "physical meaning". Indeed, we
may consider the multidimensional counterpart of the Flaschka variables
$\left(  a_{j},b_{j}\right)  $ arising from equations (\ref{ab1})-(\ref{ab3}),
which we call  \emph{generalized Flaschka variables}, \
\begin{subequations}
\begin{align}
A_{j}\left(  \theta\right)   &  :=%
{\displaystyle\sum_{k,\ell}}
\widetilde{a}_{k,\ell;j}Y_{k,\ell}\left(  \theta\right)  \qquad\text{for
}j=1,2,...,N-1\label{Aseries}\\
B_{j}\left(  \theta\right)   &  :=%
{\displaystyle\sum_{k,\ell}}
\widetilde{b}_{k,\ell;j}Y_{k,\ell}\left(  \theta\right)  \qquad\text{for
}j=1,2,...,N.\label{Aseries2}%
\end{align}
The coefficients $\widetilde{a}_{k,\ell;j}$ and $\widetilde{b}_{k,\ell;j}$ may
be found from the multidimensional Markov-Stieltjes transform, by means of the
continued fraction expansion in (\ref{CF}) for every $\left(  k,\ell\right)
,$ with $a_{j}$ and $b_{j}$ replaced by $\widetilde{a}_{k,\ell;j}$ and
$\widetilde{b}_{k,\ell;j}.$ Thus,
\end{subequations}
\[
\widetilde{a}_{k,\ell;j}=\int_{\mathbb{S}^{n-1}}A_{j}\left(  \theta\right)
Y_{k,\ell}\left(  \theta\right)  d\theta\quad\text{and }\quad\widetilde
{b}_{k,\ell;j}=\int_{\mathbb{S}^{n-1}}B_{j}\left(  \theta\right)  Y_{k,\ell
}\left(  \theta\right)  d\theta.
\]

\begin{proposition}
The series representing $A_{j}\left(  \theta\right)  $ and $B_{j}\left(
\theta\right)  $ in (\ref{Aseries})-(\ref{Aseries2}) are absolutely convergent
provided  condition (\ref{SumLambdas}) is  satisfied. 
\end{proposition}

The proof follows from Proposition \ref{PHamiltonConverge}. 

\begin{remark}
Thus, the \emph{physical interpretation} is the following: We have $N$
surfaces $A_{j},$ considered as functions on the sphere $\mathbb{S}^{n-1},$
for $j=1,2,...,N,$ which interact in the sense of classical mechanics
according to equations (\ref{ab1MULTI})-(\ref{ab3MULTI}; we see that only the
neighboring surfaces interact. 
\end{remark}

We have the following interesting representation for the Hamiltonian in
(\ref{HamiltonMulti}):
\[
H=%
{\displaystyle\sum_{k,\ell}}
H_{k,\ell}=%
{\displaystyle\int_{\mathbb{S}^{n-1}}}
{\displaystyle\sum_{j=1}^{N}}
\left(  4\times\left\vert A_{j}\left(  \theta\right)  \right\vert ^{2}%
+2\times\left\vert B_{j}\left(  \theta\right)  \right\vert ^{2}\right)
d\theta.
\]

Alternative physical meaning may be found if we use the variables $\left(
x_{k,\ell;j},y_{k,\ell;j}\right)  _{k,\ell;j}.$ 

\subsection{Physical meaning in the variables $x$ and $y$}

We have reformulated above  our pseudo-positive Toda lattice in the
generalized Flaschka  "variables" $A_{j}\left(  \theta\right)  $ and
$B_{j}\left(  \theta\right)  .$  In view of the complicated non-linear
relation between the variables $\left(  x_{j},y_{j}\right)  $ and the Flaschka
variables $\left(  a_{j},b_{j}\right)  $ it is natural to ask, if a reasonable
physical interpretation exists in the "original" variables, $x$ and $y.$ A
major criterion for success of the interpretation will be the convergence of
the series. 

Alternatively to the $\left(  a,b\right)  $ variables, we \emph{introduce} the
pseudo-positive Toda lattice by means of the generalized variables given by
the series
\begin{equation}
X_{j}\left(  x;\theta\right)  :=%
{\displaystyle\sum_{k,\ell}}
e^{-x_{k,\ell;j}}Y_{k,\ell}\left(  \theta\right)  ,\qquad\text{for
}j=1,2,...,N,\label{Xjtheta}%
\end{equation}
and
\begin{equation}
Y_{j}\left(  y;\theta\right)  :=%
{\displaystyle\sum_{k,\ell}}
y_{k,\ell;j}Y_{k,\ell}\left(  \theta\right)  \qquad\text{for }%
j=1,2,...,N.\label{Yjtheta}%
\end{equation}
We assume that the values for $x=\left(  x_{k,\ell;j}\right)  $ and $y=\left(
y_{k,\ell;j}\right)  $ are obtained from $\widetilde{a}_{k,\ell;j}$ and
$\widetilde{b}_{k,\ell;j}$ via the Flaschka change of variables in
(\ref{FlaschkaChange})-(\ref{FlaschkaChange3}). Let us recall however that the
Flaschka map is \emph{not bijective }and one has to identify equivalence
classes called configurations in the variables $\left(  x,y\right)  ,$ in
(\ref{configuration}), cf. \cite{moser1}, p. $471.$ 

The main justification for representation of the pseudo-positive  Toda lattice
in the variables $X_{j}$ and $Y_{j}$  is the following statement. 

\begin{proposition}
The series representing $X_{j}\left(  \theta\right)  $ and $Y_{j}\left(
\theta\right)  $ in (\ref{Xjtheta}) and (\ref{Yjtheta}) respectively, are
absolutely convergent in every configuration  class for $\left(  x,y\right)  $
defined by (\ref{configuration}), provided condition (\ref{SumLambdas}) is
satisfied. 
\end{proposition}%

\proof
For proving the convergence we will need explicitly the variables $x$ and $y.$
From (\ref{FlaschkaChange})-(\ref{FlaschkaChange3}) we have
\[
a_{j}=\frac{1}{2}e^{\left(  x_{j}-x_{j+1}\right)  /2},
\]
hence,
\begin{align*}
x_{1}-x_{2} &  =2\ln2a_{1}\\
x_{2}-x_{3} &  =2\ln2a_{2}\\
&  \cdot\cdot\cdot\\
x_{n-1}-x_{n} &  =2\ln2a_{n-1}.
\end{align*}
Consequently, we get the following equivalent system of equations
\begin{align*}
x_{1}-x_{2} &  =2\ln2a_{1}\\
x_{1}-x_{3} &  =2\ln2a_{1}+2\ln2a_{2}\\
&  \cdot\cdot\cdot\\
x_{1}-x_{n} &  =2\ln2a_{1}+2\ln2a_{2}+...+2\ln2a_{n-1},
\end{align*}
which implies for $j=2,3,...,n$ the expression
\[
x_{j}=x_{1}-2\ln2^{j-1}-2\ln%
{\displaystyle\prod_{m=1}^{j-1}}
a_{m}.
\]
Hence, we see that $X_{j}\left(  \theta\right)  $ defined by (\ref{Xjtheta})
satisfies
\begin{align}
X\left(  j\theta\right)   &  :=%
{\displaystyle\sum_{k,\ell}}
e^{-x_{k,\ell;j}\left(  t\right)  }Y_{k,\ell}\left(  \theta\right)
\label{XjthetaNEW}\\
&  =%
{\displaystyle\sum_{k,\ell}}
\exp\left(  -x_{k,\ell;1}+2\ln2^{j-1}+2\ln%
{\displaystyle\prod_{m=1}^{j-1}}
\widetilde{a}_{k,\ell;m}\right)  Y_{k,\ell}\left(  \theta\right)  \nonumber\\
&  =2^{2\left(  j-1\right)  }%
{\displaystyle\sum_{k,\ell}}
\left(
{\displaystyle\prod_{m=1}^{j-1}}
\widetilde{a}_{k,\ell;m}^{2}\right)  \times\exp\left(  -x_{k,\ell;1}\right)
Y_{k,\ell}\left(  \theta\right)  \nonumber
\end{align}
From the inequality between geometric and arithmetic means we obtain for every
integer $M\geq1$ the following inequality:\
\[
\left\vert
{\displaystyle\prod_{m=1}^{M}}
\widetilde{a}_{k,\ell;m}^{2}\right\vert \leq\frac{\left(  \sum_{m=1}%
^{M1}\left\vert \widetilde{a}_{k,\ell;m}\right\vert ^{2}\right)  ^{M}}{M^{M}}%
\]
From condition (\ref{SumLambdas}) we obtain
\[%
{\displaystyle\sum_{k,\ell}}
a_{k,\ell;1}^{2}<\infty.
\]
On the other hand, recall that the Flaschka map (\ref{FlaschkaChange}%
)-(\ref{FlaschkaChange3}) is not on-to-one and all $x_{k,\ell;j}$ are only
determined up to a constant $C_{k,\ell},$ i.e. $\left(  x_{k,\ell;j}%
+C_{k,\ell}\right)  _{j=1}^{N}$ belong to the same class; hence, we may
choose $x_{k,\ell;1}$ in such a way that
\[
\left\vert e^{-x_{k,\ell;1}}Y_{k,\ell}\left(  \theta\right)  \right\vert \leq
C\qquad\text{for all }\left(  k,\ell\right)  ;
\]
in particular we may choose $e^{-x_{k,\ell;1}}=k^{-\left(  n-2\right)  },$ or
$x_{k,\ell;1}=\ln\left(  k^{n-2}\right)  .$ Thus for every class of
equivalence we find a convergent series  for $X_{j}\left(  \theta\right)  .$ 

The convergence of the series for $Y_{j}\left(  \theta\right)  $ follows
easily, since $y_{k,\ell;j}$ are proportional to the variables $\widetilde
{b}_{k,\ell;j}$ which are bounded by (\ref{Htilde}) and condition
(\ref{SumLambdas}).%

\endproof

Finally, we comment on a physical interpretation of the above representation. 

\begin{remark}
In analogy to the classical Toda lattice, we may assign some \textbf{physical
meaning }(although not so directly understood)\textbf{ }to the set of surfaces
$X_{j}\left(  \theta\right)  $ and $Y_{j}\left(  \theta\right)  $ over the
sphere, by considering them as mass-surfaces which interact by means of the
Hamiltonian system for the infinite set of variables $\left(  x_{k,\ell
;j},y_{k,\ell;j}\right)  _{k,\ell;j}$ provided in (\ref{TodaLattice}), with
Hamiltonian (\ref{HamiltonMulti}). Let us remark that these equations may be
considered as a non-linear system of pseudodifferential equations on the
sphere with respect to the functions $X\left(  j\theta\right)  =X_{j}\left(
\theta\right)  $ and $Y\left(  j\theta\right)  =Y_{j}\left(  \theta\right)  .$ 
\end{remark}

\section{A class of isospectral measures \label{Sspecial}}

There is an important class of measures which has been identified in
\cite{kounchevrender}:

\begin{definition}
\label{Dintegrability}We say that the pseudo-positive measure $\mu$ satisfies
the integrability condition iff the following inequality holds
\begin{equation}
\sum_{k,\ell}\int_{0}^{\infty}\frac{d\mu_{k,\ell}\left(  r\right)  }{r^{k}%
}<\infty\qquad\text{for }k\in\mathbb{N}_{0},\ \ell=1,...,d_{k}.
\label{integrabilityCondition}%
\end{equation}

\end{definition}

In \cite{kounchevrender} it has been proved that for such measures we may
construct a Gauss-Jacobi measure on $\mathbb{R}^{n}$ which provides a
generalization of the Gauss quadrature formula on $\mathbb{R}$. The new
formula has been called polyharmonic Gauss-Jacobi cubature formula.

Further we will consider measures $\mu\left(  x;t\right)  $ which depend on
both the space variable $x$ and the time variable $t$ and which have the
isospectral property in every component $\mu_{k,\ell}\left(  r;t\right)  ,$
namely for every index $\left(  k,\ell\right)  $ we consider the measure on
$\mathbb{R}_{+}$
\[
d\mu_{k,\ell}\left(  \rho;t\right)  =\sum_{j=1}^{N}r_{k,\ell,j}^{2}\left(
t\right)  \delta\left(  \rho-\lambda_{k,\ell,j}\left(  t\right)  \right)
,\qquad\rho\in\mathbb{R}_{+},t\in\mathbb{R},
\]
where the functions $r_{k,\ell,j}\left(  t\right)  \geq0$ and $\lambda
_{k,\ell,j}\left(  t\right)  \geq0$ satisfy the equations of the Toda flow:\
\[
\frac{d}{dt}\lambda_{k,\ell,j}\left(  t\right)  =0,\qquad\frac{d}{dt}%
r_{k,\ell,j}\left(  t\right)  =-\lambda_{k,\ell,j}\left(  t\right)
r_{k,\ell,j}^{2}\left(  t\right)  .
\]
We see from the first equation that $\lambda_{k,\ell,j}$ are constants (with
respect to the variable $t$), hence the second equation becomes
\begin{equation}
\frac{d}{dt}r_{k,\ell,j}\left(  t\right)  =-\lambda_{k,\ell,j}r_{k,\ell,j}%
^{2}\left(  t\right)  \label{EqMotion}%
\end{equation}
with $\lambda_{k,\ell;j}=\lambda_{k,\ell;j}\left(  0\right)  $ for
$t\in\mathbb{R}.$ We consider the measure defined by the formal series
\[
d\mu\left(  x;t\right)  =%
{\displaystyle\sum_{k,\ell}}
d\mu_{k,\ell}\left(  \rho;t\right)  Y_{k,\ell}\left(  \theta\right)
\]
for which it is not \emph{a priori} clear that the series is convergent in any
sense. We prove its convergence in the following theorem.

\begin{theorem}
\label{Tintegrability} Let the measure $d\mu\left(  x;0\right)  $ satisfy the
integrability condition (\ref{integrabilityCondition}). Then for every $t>0$
the measure $d\mu\left(  x;t\right)  $ satisfies the integrability condition
(\ref{integrabilityCondition}).
\end{theorem}

%

\proof
Let us consider
\[
S_{k,\ell}\left(  t\right)  :=\int_{0}^{\infty}\frac{d\mu_{k,\ell}\left(
\rho;t\right)  }{\rho^{k}}=%
{\displaystyle\sum_{j=1}^{N}}
\frac{r_{k,\ell,j}^{2}\left(  t\right)  }{\lambda_{k,\ell,j}^{k}},\qquad
t\in\mathbb{R}.
\]
We obtain for its derivative, using the equation of motion (\ref{EqMotion}),
\begin{align*}
\frac{d}{dt}S_{k,\ell}\left(  t\right)   &  =\frac{d}{dt}%
{\displaystyle\sum_{j=1}^{N}}
\frac{r_{k,\ell,j}^{2}\left(  t\right)  }{\lambda_{k,\ell,j}^{k}}=2%
{\displaystyle\sum_{j=1}^{N}}
\frac{r_{k,\ell,j}\left(  t\right)  }{\lambda_{k,\ell,j}^{k}}\frac{d}%
{dt}r_{k,\ell,j}\left(  t\right)  =-2%
{\displaystyle\sum_{j=1}^{N}}
\frac{r_{k,\ell,j}\left(  t\right)  }{\lambda_{k,\ell,j}^{k}}\lambda
_{k,\ell,j}r_{k,\ell,j}^{2}\left(  t\right) \\
&  =-2%
{\displaystyle\sum_{j=1}^{N}}
\frac{r_{k,\ell,j}^{3}\left(  t\right)  }{\lambda_{k,\ell,j}^{k-1}}%
,\qquad\text{for }k\in\mathbb{N}_{0},\ \ell=1,...,d_{k}.
\end{align*}
By definition the functions $r_{k,\ell,j}\left(  t\right)  $ and
$\lambda_{k,\ell,j}$ are non-negative, hence it follows that $S_{k,\ell
}\left(  t\right)  $ is a decreasing function of $t.$ Hence,
\[%
{\displaystyle\sum_{k,\ell}}
\int_{0}^{\infty}\frac{d\mu_{k,\ell}\left(  \rho;t\right)  }{\rho^{k}}%
\]
is a decreasing function of $t$ which proves our statement, according to
Definition \ref{Dintegrability}.%

\endproof

Acknowledgements. The hospitality of IAM and HCM of the University of Bonn is
gratefully acknowledged, as well as the  support by the Alexander von Humboldt Foundation.

\end{document}